\documentclass{article}

\usepackage{PRIMEarxiv}

\usepackage[utf8]{inputenc} % allow utf-8 input
\usepackage[T1]{fontenc}    % use 8-bit T1 fonts
\usepackage{hyperref}       % hyperlinks
\usepackage{url}            % simple URL typesetting
\usepackage{booktabs}       % professional-quality tables
\usepackage{amsmath}
\usepackage{amsthm}
\usepackage{amsfonts}       % blackboard math symbols
\usepackage{nicefrac}       % compact symbols for 1/2, etc.
\usepackage{microtype}      % microtypography
\usepackage{tabularx}
\usepackage{multirow}
\usepackage[ruled,linesnumbered]{algorithm2e}
\usepackage{fancyhdr}       % header
\usepackage{graphicx}       % graphics
\graphicspath{{media/}}     % organize your images and other figures under media/ folder

\newcolumntype{Y}{>{\centering\arraybackslash}X}

\newtheorem{remark}{Remark}[section]

\title{AutoAMG($\theta$): An Auto-tuned AMG Method Based on Deep Learning for Strong Threshold
}

\author{
  Haifeng Zou \\
  Graduate School of China Academy of Engineering Physics, China Academy of Engineering Physics \\
  \texttt{zouhaifeng19@gscaep.ac.cn} \\
  \And
  Xiaowen Xu \\
  Laboratory of Computational Physics, Institute of Applied Physics and Computational Mathematics \\
  \texttt{xwxu@iapcm.ac.cn} \\
  \And
  Chen-Song Zhang \\
  Academy of Mathematics and Systems Science, Chinese Academy of Sciences \\
  \texttt{zhangcs@lsec.cc.ac.cn} \\
  \And
  Zeyao Mo \\
  Laboratory of Computational Physics, Institute of Applied Physics and Computational Mathematics \\
  \texttt{xwxu@iapcm.ac.cn} \\
}

\begin{document}
\maketitle

\begin{abstract}
Algebraic Multigrid (AMG) is one of the most used iterative algorithms for solving large sparse linear equations $Ax=b$. In AMG, the coarse grid is a key component that affects the efficiency of the algorithm, the  construction of which relies on the strong threshold parameter $\theta$. This parameter is generally chosen empirically, with a default value   in many current AMG solvers of 0.25 for 2D problems and 0.5 for 3D problems. However, for many practical problems, the quality of the coarse grid and the efficiency of the AMG algorithm are sensitive to $\theta$; the default value is rarely optimal, and sometimes is far from it. Therefore, how to choose a better $\theta$ is an important question. In this paper, we propose a deep learning based auto-tuning method, AutoAMG($\theta$) for multiscale sparse linear equations, which are widely used in practical problems. The method uses Graph Neural Networks (GNNs) to extract matrix features, and a Multilayer Perceptron (MLP) to build the mapping between matrix features and the optimal $\theta$, which can adaptively output $\theta$ values for different matrices. Numerical experiments show that AutoAMG($\theta$) can achieve significant speedup compared to the default $\theta$ value.
\end{abstract}

% keywords can be removed

\section{Introduction}
\label{sec1}

Solving sparse linear equations $Ax=b$ is ubiquitous in numerical simulations, and is   a major bottleneck affecting computational efficiency. Owing  to its good generality and optimal computational complexity, the AMG algorithm~\cite{ruge1987algebraic,stuben2001review,xu2017} is one of the most widely used algorithms for large-scale sparse linear equations, which uses only information from the matrix to construct components, including coarsening, interpolation, and restriction operators. During the coarsening procedure, a subset of points from the adjacency matrix $A$ is selected as points in the coarse grid, which is the basis for constructing a coarse grid matrix $A_c$. Different coarsening strategies will result in different coarse matrices $A_c$. In the classical AMG algorithm, points in the subset are selected based on the strong threshold $\theta$ and the strength of the connectivity between points, which is calculated by the value of the matrix entries. Hence the value of $\theta$ directly affects the grid coarsening results, and is a key factor affecting the algorithm's  efficiency.

In the classical AMG algorithm, most  coarsening algorithms are based on heuristic strategies for coarse grid construction. A basic principle is to perform coarsening along the direction of strong connectivity  to accommodate the property that algebraic errors are smoothed or relaxed along the same direction. If the strong threshold $\theta$ is large, then the number of points in the corresponding coarse grid is large, which means the AMG algorithm has high complexity. If $\theta$ is small, although the number of points in the coarse grid is smaller, the residuals may decrease more slowly, requiring more iterations to converge. Since there is no strict theoretical guarantee on the size of the optimal coarse grid, the current value of $\theta$ can only be chosen empirically. For example, in the  HYPRE AMG solver~\cite{falgout2002hypre}, depending on the physical dimension of the sparse matrix, $\theta$ equals 0.25 for 2D problems and 0.5 for 3D problems. However Vakili~\cite{vakili2009recommendations} and Nikola~\cite{kosturski2015performance} utilize the incompressible Navier Stokes equation and linear poroelasticity equation, respectively, as the test cases, both of which  show  the increase of $\theta$ along with the monotone decrease of time. Here, we take the diffusion problem as the example, and find that the number of iterations changes irregularly  with the increase of $\theta$. If the diffusion coefficients are isotropic, the default values of $\theta$ can achieve the desired convergence rate. If the diffusion coefficients are anisotropic, which means there are significant differences in the strength of connectivity between points, then the default values of $\theta$ maybe far from the optimal. Notably, small changes in $\theta$ may have a large impact on the construction of the coarse grid, thus affecting the convergence rate and efficiency of AMG. In particular, we focus on the so-called multiscale sparse matrices~\cite{xu2017algebraic}. In some typical test cases, the number of iterations of the default $\theta$ is 10 times larger than the minimum number of iterations obtained by grid search (see Section \ref{sec2-3} for detail).

The above problem can be summarized as follows: how to choose an appropriate $\theta$ for any given sparse matrix. Considering that the properties of the input matrix may vary dynamically, the automatic selection of a suitable $\theta$ for different linear systems is a crucial and challenging task, since there is no theoretical guarantee yet. Machine learning and deep learning algorithms provide a feasible approach. Paola F~\cite{Antonietti2021} used a Convolutional Neural Network (CNN) to extract matrix features, and built a regression model with those features. The inputs of the regression model are matrix features, strong threshold $\theta$, and $-\log_{2} h$ ($h$ is the edge length in the mesh), and the output $y$ is an approximated convergence factor. After training, the regression model is used to optimize $\theta$. There are other ways to enhance the robustness of iterative methods with machine learning and deep learning. For example, a variety of classification algorithms are used to select optimal iterative methods based on the input matrix features~\cite{Holloway2007,bhowmick2006application,Eller2012,motter2015lighthouse}; deep learning algorithms are utilized to optimize the prolongation matrix $P$, restriction matrix $R$, and smoother $S$ in AMG~\cite{katrutsa2017deep, Greenfeld2019, Luz2020, chen2020meta}.

Our target is optimizing $\theta$ adaptively according to the input matrices, and our contributions are as follows:
\begin{itemize}
	\item Classical graph convolution networks such as GCN~\cite{kipf2017semisupervised}, GIN~\cite{xu2018powerful} are used to extract matrix features, but they didn't work well. Therefore, a new variant of graph convolutional network is proposed in this paper as the feature extractor (see Section \ref{sec3-gcin} for detail).
	\item We utilize MLP to directly build the mapping between matrix features and the optimal $\theta$, avoiding optimizing the regression model.
\end{itemize}
The strong threshold $\theta$ auto-tuning method is called AutoAMG($\theta$), and its effectiveness is verified by matrices from the diffusion equations and radiation diffusion equations~\cite{xiaowen2009algebraic,xu2017algebraic}. Numerical experiments show that AutoAMG($\theta$) can achieve acceleration by a factor of 4.47 relative to the default $\theta$ in diffusion equations, and a factor of 11.63 relative to the default $\theta$ in radiation diffusion equations.

The rest of this paper is organized as follows. Section \ref{sec2} briefly introduces the rationale behind AMG and shows how $\theta$ affects iteration. Section \ref{sec3} explains the details of AutoAMG($\theta$). Section \ref{sec4} presents numerical experiments and results of AutoAMG($\theta$). Section \ref{sec5} summarizes our work.

\section{Sensitivity of strong threshold}
\label{sec2}

\subsection{AMG algorithm}
\label{sec2-1}

The AMG algorithm can be divided into two phases: SETUP and SOLVE, as described in Algorithms \ref{alg-setup} and \ref{alg-tg}, respectively. Considering the complexity of AMG, we introduce its simplified version, the Two-Grid~(TG) algorithm.

In the SETUP phase, the TG algorithm constructs a coarse-level grid, an interpolation matrix $P$, and a restriction matrix $R$ based on the matrix $A$. In the SOLVE phase, it performs a standard multigrid cycle based on the matrices generated in SETUP phase, including pre-smoothing, restricting residuals to the coarse grid, solving residual equations in the coarse grid, interpolating the error back to the fine-level grid for correction, and post-smoothing. In particular, if the TG algorithm is called recursively to solve linear equations in the coarse-level grid (line 6, Algorithm \ref{alg-tg}), it becomes a multigrid algorithm.

\begin{algorithm}
	Coarsening: Construct the fine-level grid based on the matrix $A$ and let $\Omega$ be the set containing all fine-level variables. Split the set $\Omega$ into  set $C$ containing all coarse-level variables and  set $F$ containing the remaining fine-level variables, according to the strong threshold $\theta$. In addition, $F \cap C= \emptyset$, $F \cup C = \Omega$. \\
	
	Computing $A_c$: Based on the coarse variable set $C$, compute the interpolation matrix $P$ and   restriction matrix $R$. Then compute the coarse-level matrix $A_c$ by $A_{c}=RAP$. \\
	
	\caption{SETUP phase}
	\label{alg-setup}
\end{algorithm}

\begin{algorithm}[H]
	Pre-smoothing: smoothing $\mu_1$ times on $Ax=b$, get the approximate solution $x_f$ \\
	\eIf{deepest level }{
		Solve $Ax=b$ directly	
	} 
	{
		Restricting residuals into coarse grid: $b_c = R(b-Ax_f)$\\
		
		Solving the coarse grid equation: $A_c x_c = b_c$\\
		
		Interpolating and correcting: $x_f = x_f + Px_c$\\
	}
	Post-smoothing: smoothing $\mu_2$ times on $Ax=b$, update $x_f$
	\caption{SOLVE phase}
	\label{alg-tg}
\end{algorithm}
\mbox{}

In the SETUP phase of the classical AMG, the algorithm will split all variables into a coarse variable set $C$ and fine variable set $F$ (C/F splitting), which is the first step in   Algorithm \ref{alg-setup}. More specifically, let $N_i = \{j \ | \ a_{ij} \neq 0, j \neq i \}$ be the dependency set of variable $i$, i.e.,  $i$ strongly depends on  $j$ (or  $j$ strongly influences  $i$). If
\begin{align}
|a_{ij}| \geq \theta \max_{k\in N_i,k\neq i} |a_{ik}|, \label{eq-theta}
\end{align}
where $0<\theta \leq 1$ is the strong threshold, then we can define the strong dependency set $S_i$ and   strong influence set $S_i^T$ of variable $i$,
\begin{align*}
S_i &= \left\lbrace j \ \bigg| \ |a_{ij}| \geq \theta \max_{k\in N_i,k\neq i} |a_{ik}|, \  j \in N_i  \right\rbrace  , \\
S_i^T &= \left\lbrace j \ \bigg| i \in S_j, \  j \in N_i  \right\rbrace .
\end{align*}

According to the definitions, a basic principle of coarsening is that the larger $|S_i^T|$ is, the more important the variable $i$ is, and the more likely it is to be selected as a coarse variable. Following this principle, the result of grid coarsening is closely related to the strength of connectivity between variables, i.e., it relies on the strong threshold $\theta$ in Eq. (\ref{eq-theta}).

\subsection{Multiscale matrix}
\label{sec2-2}

Multiscale matrices are common in practical problems. Factors such as multimedia (e.g., anisotropy, discontinuity, oscillating coefficients), large deformations, strong nonlinearities, and multiphysics coupling all lead to the multiscale property of   matrices obtained by discretization. Define the matrix $A \in \mathbb{R}^{n \times n}$, and let $\Omega = \{0,1,2,\cdots,n \}$ be the set containing all row indices of the matrix. Given a multiscale threshold $\delta \geq 0$, define the   multiscale set  
\begin{align}
\Omega_{MS} = \left\lbrace i \ \bigg| \ i \in \Omega, \ \log_{10}(\frac{\max_{k\in N_i,k\neq i} |a_{ik}|}{\min_{k\in N_i,k\neq i} |a_{ik}|} ) \geq \delta  \right\rbrace  . \label{eq-ms}
\end{align}
If $\Omega_{MS} \neq \emptyset$, then $A$ is defined as a multiscale matrix (under the threshold $\delta$). If $\Omega_{MS} = \emptyset$, then $A$ is a single-scale matrix. 

A detailed definition of the multiscale matrix and how the multiscale property affects the AMG algorithm can be found in~\cite{xu2017algebraic}. From Eq. (\ref{eq-ms}), the multiscale property reflects the strength of the numerical difference between the maximum and minimum absolute values of the nondiagonal elements in the same row of the matrix.

\subsection{Impact of $\theta$}
\label{sec2-3}

The effect of $\theta$ on the efficiency of the AMG algorithm is illustrated by the diffusion equation below,
\begin{equation}
\begin{split}
-  \nabla \cdot (\kappa \nabla  u) &= f_1, \quad x \in \Omega \, ,\\
u &= f_2, \quad x \in \partial \Omega \, , \\
\end{split}
\label{eq-1}
\end{equation}
where $\kappa$ is the diffusion coefficient. In a two-dimensional (2D) diffusion problem, we define the diffusion coefficient as 
\begin{align*}
\kappa = \begin{bmatrix}
10^{4\varepsilon} & 0 \\
0 & 1 \\
\end{bmatrix} \, ,
\end{align*}
where $0 \leq \varepsilon \leq 1$ is a random number. 

\begin{figure}[h]
	\centering
	\includegraphics[scale=0.5]{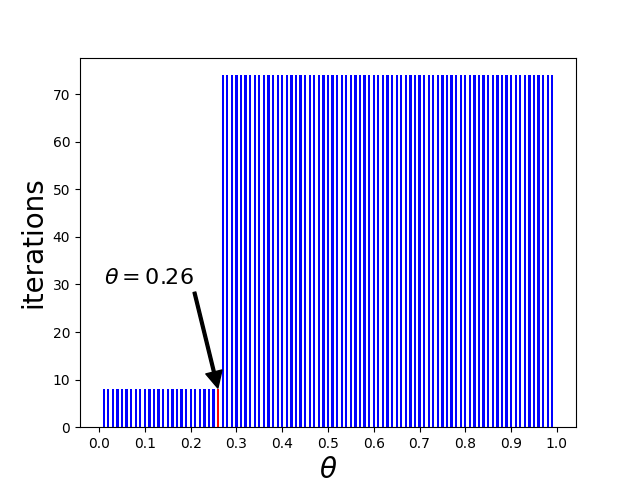}
	\caption{"Boundary" matrix with 144 rows, 218 nonzeros}
	\label{fig-sen}
\end{figure}

We theoretically verify the effect of $\theta$  based on a specific small matrix whose inverse we can compute. The matrix comes from diffusion equation (\ref{eq-1}), with a random diffusion coefficient $\kappa$ and a mesh size of $12\times 12$. We gradually drop the element with the minimum absolute value in the matrix to obtain a "boundary" matrix that is one step from a single-scale matrix. Then the matrix is solved by the TG algorithm, with results as shown in Figure \ref{fig-sen}, where the $x$-axis is $\theta$, in the interval $(0,1)$, with a common difference of 0.01, and the $y$-axis is the number of iterations (upper limit   500). This shows that there is a critical value $\theta^*=0.26$ in the matrix, where the number of iterations is 8 when $\theta \leq \theta^*$,  and 74 when $\theta > \theta^*$.

Based on the analysis of the convergence factor in the TG algorithm~\cite{falgout2005two}, we compare the theoretically estimated and computed convergence factors in Table \ref{tab-sen}. The theoretical results remain consistent with the computed results, which indicates that the phenomenon of an oscillating number of iterations is caused by the algorithm itself, and is an essential feature of the algorithm. Such results further illustrate the necessity of optimizing $\theta$.
\begin{table}[h!]
	\centering
	\caption{Theoretical and computed convergence factors}
	\label{tab-sen}
	\begin{tabularx}{0.9\textwidth}{|Y|Y|Y|}
		\hline 
		$\theta$ & Theoretical & Computed \\
		\hline 
		0.26	& 0.2500  & 0.2498 \\
		\hline 
		0.27  & 0.9477  & 0.9477 \\
		\hline 
	\end{tabularx}
\end{table}

\begin{figure}
	\begin{minipage}[t]{0.5\linewidth}
		\centering
		\includegraphics[width=2.6in]{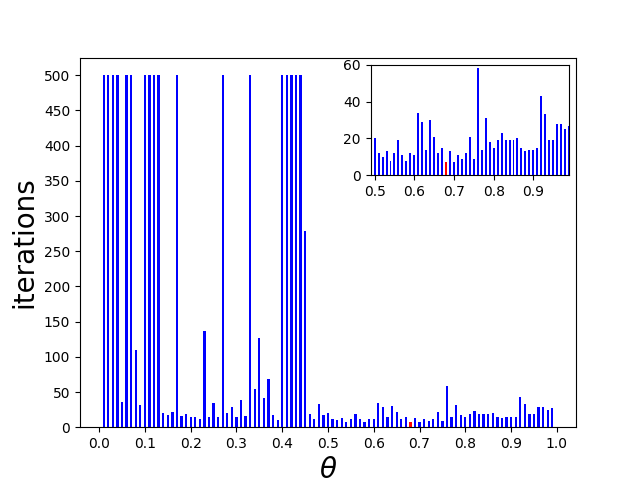}
	\end{minipage}%
	\begin{minipage}[t]{0.5\linewidth}
		\centering
		\includegraphics[width=2.6in]{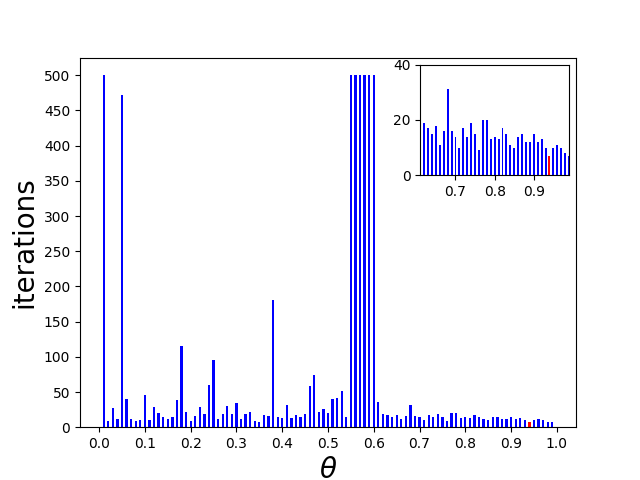}
	\end{minipage}
	\caption{Two random seeds with the same DoF = 1048576}
	\label{fig-2d}
\end{figure}

Furthermore, the numbers of iterations based on two   random seeds are depicted in Figure~\ref{fig-2d}, where the mesh size is $1024\times 1024$, with 1048576 degrees of freedom (DoF). The iterative method is GMRES, with AMG as the precondition and PMIS~\cite{luby1985simple} as the coarsening algorithm. Figure \ref{fig-2d} shows that first, for both random coefficients, the number of iterations changes irregularly with $\theta$, and second, different random seeds have different behaviors.

Table \ref{tab-2d} shows the maximum and minimum number of iterations for these cases, as well as the number of iterations corresponding to the default $\theta$. The maximum and minimum number of iterations for both cases are 500 and 7, which means there is a large gap between the maximum and minimum. Moreover, the values of $\theta$ corresponding to the maximum are not the same (0.68 and 0.94). Concerning the default value $\theta=0.25$, the number of iterations is 35 and 95 in two cases, which are 5 and 13 times larger than the corresponding minimum. These results also imply that for random diffusion coefficients, the number of iterations is sensitive to the value of $\theta$, and the value of $\theta$ corresponding to the minimum is different for different matrices.

\begin{table}[h]
	\centering
	\caption{Min, Max, Default iterations and corresponding $\theta$}
	\label{tab-2d}
	
	\begin{tabularx}{0.9\textwidth}{|Y|Y|Y|Y|}
		\hline
		& Min/$\theta$ & Max/$\theta$ & Default ($\theta = 0.25$) \\
		\hline 
		Left &7 / 0.68 & 500 / 0.01 & 35\\
		\hline 
		Right & 7 / 0.94& 500 / 0.01& 95 \\
		\hline 
	\end{tabularx}
\end{table}

\section{AutoAMG($\theta$): auto-tune $\theta$ for multiscale matrices}
\label{sec3}

\subsection{AutoAMG($\theta$) procedure}
\label{sec3-1}

\begin{figure}[h] 
	\centering
	\includegraphics[width=0.9\linewidth]{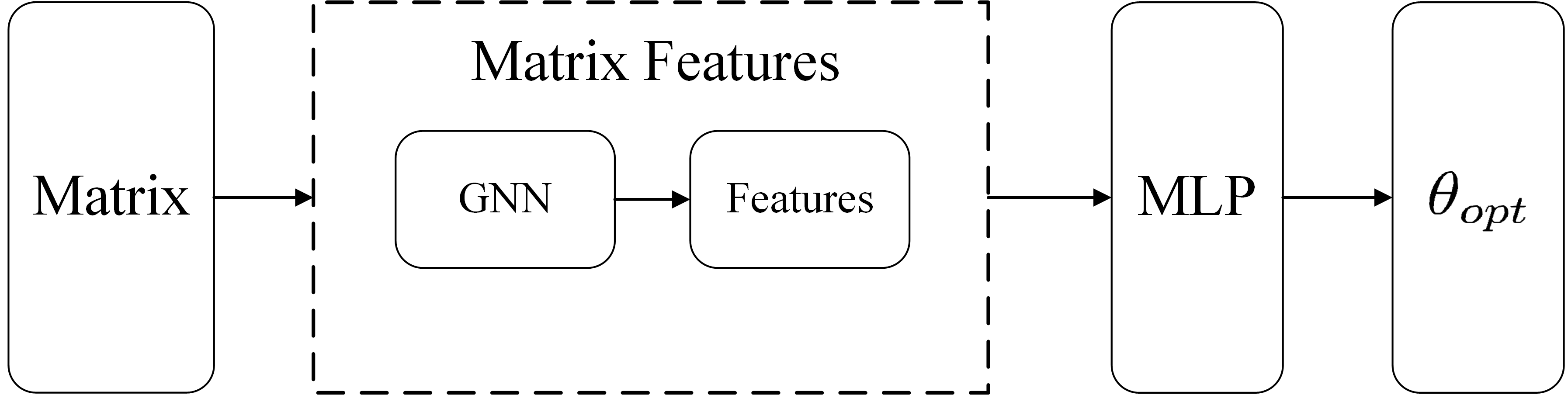}
	\caption{AutoAMG($\theta$) procedure}
	\label{fig-auto}
\end{figure}

The comprehensive AutoAMG($\theta$) procedure is depicted in Figure \ref{fig-auto}. The input of AutoAMG($\theta$) is the matrix, which is treated as the adjacent graph. GNN based on message passing is utilized to extract graph features. Subsequently, AutoAMG($\theta$) establishes a mapping between these extracted features and the optimal value of $\theta_{opt}$. Note that $\theta_{opt}$ pertains to the $\theta$ value yielding the fewest iterations during grid search.

The key step in AutoAMG($\theta$) is feature extraction. Considering matrices discretized from the same equation, their sparsity patterns exhibit a degree of similarity, differing in the number of rows and element values. Notably, the comparison depicted in Figure \ref{fig-2d} illustrates that conventional structural and numerical matrix features (e.g., dimensions, sparsity patterns) fall short of adequately capturing the intricate influence of $\theta$ on the iterations across diverse matrices. Besides, the calculation of spectral attributes (e.g., condition number, eigenvalue distribution) is time-consuming, sometimes even surpassing the time required for solving the linear equation. In AutoAMG($\theta$),  GNN is utilized to extract node features in graphs, then graph features are derived based on the extracted node features.

\subsection{GNN}
GNN is one of the deep learning algorithms specifically designed for the analysis of graph data structure. Nowadays, GNNs are utilized in diverse domains such as social recommendation, traffic prediction, and molecular structure prediction, et al~\cite{wu2020comprehensive}. A graph $G$ is represented as $G=(V, E)$, where $V$ is the set of nodes and $E$ is the set of edges in the graph. The number of nodes is $|V|$ and the number of edges is $|E|$. Let $v_i \in V$ denote the $i$-th node and $e_{ij}=(v_i,v_j) \in E$ denote the directed edge from node $v_i$ to node $v_j$. Let $N(v)$ denote all neighbor nodes of the node $v$. Since every node and edge may have features, let $X_v \in \mathbb{R}^{|V|\times d}$ denote feature matrix of all nodes and $X_e \in \mathbb{R}^{|E|\times c}$ denote feature matrix of all edges, where $X_{v_i} \in \mathbb{R}^{d}$ is the feature vector of the $i$-th node and $X_{e_{ij}} \in \mathbb{R}^{c}$ is the feature vector of edge $e_{ij}$. Let $X_g$ denote the feature vector of the graph.

The standard operation of a GNN involves the following process: commencing with the initial node feature vector $X_v^{(0)}$ and edge feature vector $X_e^{(0)}$, diverse GNN variants employ distinct strategies to iteratively update these feature vectors for nodes and edges. This evolution is often visualized as a mechanism of message passing that transpires among the nodes within the graph, whose formula is\footnote{https://pytorch-geometric.readthedocs.io/en/latest/tutorial/create\_gnn.html}
\begin{align}
X_{v_i}^{(k)} = \gamma^{(k)} \left( X_{v_i}^{(k-1)}, Aggr^{(k)}_{j\in N(i)}  \ \phi^{(k)}(X_{v_i}^{(k-1)},X_{v_j}^{(k-1)},X_{e_{ji}}) \right) \ , \label{eq-mp}
\end{align}
where $\phi^{(k)}$, $Aggr$ and $\gamma^{(k)}$ are three kernel functions of the GNN algorithm:
\begin{itemize}
	\item $\phi^{(k)}$ is the message function that dictates the content of messages propagated by the neighboring nodes and edges of node $v_i$;  
	\item $Aggr^{(k)}$ is the aggregation function that defines the approach taken to process the sent messages;
	\item $\gamma^{(k)}$ is the update function that specifies how the node feature vector $X_{v_i}^{(k-1)}$ and the aggregated messages are combined to derive the updated node feature vector $X_{v_i}^{(k)}$.
\end{itemize} 
These three functions can either be differentiable functions or MLPs. Each message passing step corresponds to a GNN layer, and these functions may vary across different layers. After $K$ steps, the resultant node feature vector $X_{v_i}^{(K)}$ is used for downstream tasks, such as node classification. It's worth noting that Eq. (\ref{eq-mp}) focuses on the nodes within the graph, while there exist GNNs that involve the updating of the edge feature vector $X_{e_{ji}}$~\cite{wang2019dynamic}. Utilizing $X_{v_i}^{(K)}$, the computation of the graph feature vector $X_g$ is facilitated via a Readout function. A variety of Readout functions are available for selection, such as the SUM function
\begin{align}
X_g = \sum_{i=1}^{|V|} X_{v_i}^{(K)} \ , \label{eq-read1}
\end{align}
which is the sum of all nodes features; or MEAN function
\begin{align}
X_g = \frac{1}{|V|} \sum_{i=1}^{|V|} X_{v_i}^{(K)} \ , \label{eq-read2}
\end{align}
which is the average of all nodes features, .et al.

At first, we tried to use GCN~\cite{kipf2017semisupervised} and GIN~\cite{xu2018powerful} to extract graph features. According to Eq. (\ref{eq-mp}), a single GCN layer is defined as 
\begin{align}
X_{v_i}^{(k)} = MLP^{(k)} \left(  \sum_{v_j \in N(v_i) \cup v_i } \frac{w_{ji}}{\sqrt{\hat{D}_i \hat{D}_j}} X_{v_j}^{(k-1)}  \right)  \ , \label{eq-gcn}
\end{align}
where $w_{ji}$ is the weight of edge $e_{ji}$, if the graph is unweighted, then $w_{ji}=1$; $\hat{D}$ is the diagonal degree matrix and $\hat{D}_i$ is the degree of node $v_i$ in the graph. In GCN, $\phi^{(k)}$ is the feature vector of neighbor nodes, $Aggr^{(k)}$ is weighted average, and $\gamma^{(k)} = MLP^{(k)}$. MEAN function (Eq. (\ref{eq-read2})) commonly serves as the Readout function of GCN.

A single GIN layer is defined as 
\begin{align}
X_{v_i}^{(k)} = MLP^{(k)} \left( w_{ii} (1+\epsilon^{(k)}) \cdot X_{v_i}^{(k-1)} + \sum_{v_j \in N(v_i)} w_{ji} X_{v_j}^{(k-1)}  \right)  \ , \label{eq-gin}
\end{align}
where $w_{ii}$ is the weight of node $v_i$'s self loop, $w_{ji}$ is the weight of edge $e_{ji}$, and $\epsilon^{(k)}$ can be a trainable parameter or a fixed constant number. Compared to Eq. (\ref{eq-gcn}), the aggregation function $Aggr^{(k)}$ is summation. The authors of GIN demonstrated that in specific scenarios, the MEAN and MAX functions would impair the expressiveness of the GNN. Consequently, both the aggregation and Readout functions in GIN are summation rather than average. The recommended Readout function for GIN is
\begin{equation}
\begin{split}
X_g = \text{CONCAT} &\left( \text{SUM} \left( X_v^{(k)} \right)   \bigg| k=0,1,\dots K \right) \ , \\
\text{SUM} \left( X_v^{(k)} \right)  &=  \sum_{i=1}^{|V|} X_{v_i}^{(k)} \ ,
\end{split} \label{eq-gin2}
\end{equation}
where CONCAT is the concatenation function that concatenate several vectors into a long vector.

\subsection{GCIN}
\label{sec3-gcin}

We choose GCN and GIN from the existing GNNs for matrix feature extraction due to their low computational complexity ($\mathcal{O}(N)$). Moreover, each layer can be implemented using Sparse Matrix-Vector Multiplication (SpMV) operations, facilitating integration of these GNNs into existing iterative software frameworks.

However, our experimental results revealed that GCN and GIN did not yield satisfactory outcomes. The issue with GCN was the occurrence of NAN (Not A Number) errors during the training phase. Upon conducting a thorough debugging process, we identified the source of these NAN errors to be the degree matrix $\hat{D}$ in Eq. (\ref{eq-gcn}). Specifically, the degrees of certain nodes were either NAN or INF, thereby leading to $\hat{D}_i^{-1/2} = NAN$. The GCN implementation was based on PyTorch Geometric~\cite{2019pyg}. Despite our attempts to use the latest software version, the encountered errors persisted unchanged.

\begin{figure}[h] 
	\centering
	\includegraphics[width=0.5\linewidth]{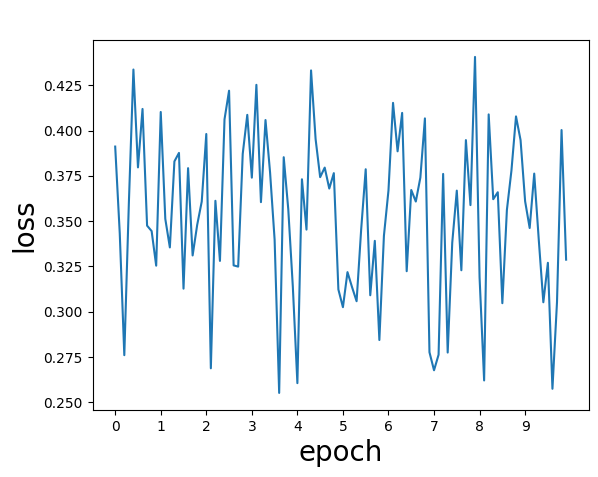}
	\caption{The training process of GIN}
	\label{fig-gin}
\end{figure}

The problem encountered with GIN pertained to the absence of a reduction in the loss value during training, as shown in Figure \ref{fig-gin}. This phenomenon is plausible given the nature of this problem, where the absence of normalization in GIN (refer to Eq. (\ref{eq-gin}) and Eq. (\ref{eq-gin2})) allows values to accumulate, consequently impeding the convergence process.

The experiments of GIN reveal that normalization is essential for our problem. Nonetheless, improper normalization can lead to NAN errors during training. After testing and analyzing, we introduce the Graph Convolutional Isomorphism Network (GCIN), which amalgamates the attributes of both GCN and GIN. A single layer of GCIN is defined as 
\begin{align}
X_{v_i}^{(k)} = MLP^{(k)} \left(   \sum_{v_j \in N(v_i)\cup v_i } w_{ji} X_{v_j}^{(k-1)}  \right) \ , \label{eq-gcin1}
\end{align}
and the Readout function is 
\begin{align}
X_g = \sum_{k=1}^{K} \left( \frac{1}{|V|} \sum_{i=1}^{|V|} X_{v_i}^{(k)} \right) \ . \label{eq-gcin}
\end{align}
Notably, normalization is integrated within the Readout function rather than being incorporated into the message passing process.

\subsection{Optimizing strong threshold $\theta$}
\label{sec3-3}

Following the extraction of matrix features, the subsequent phase involves the optimization of the strong threshold $\theta$. A conventional approach encompasses training a regression model, where matrix features and $\theta$ are inputs, and the performance metric (such as computation time, iteration count, or convergence factor) serves as the output. Then the optimization of $\theta$ relies on this regression model. Here, let the graph feature vector $X_g$ denote the matrix features, $y$ denote the performance metric, and $f$ denote the regression function. Consequently, the regression model is expressed as follow 
\begin{align}
y = f(X_g,\ \theta). \label{eq-reg}
\end{align}
Upon completion of the training phase, the regression function $f$ is established. Given any matrix, the optimization problem can be written as
\begin{equation*}
\begin{split}
\max_{\theta\in (0,1)} \ &y = f(X_g,\ \theta) \ ,
\end{split} 
\end{equation*}
which is a black-box optimization problem. To circumvent the need for solving this problem, we forego the creation of a regression model like Eq. (\ref{eq-reg}), opting to establish a direct mapping between matrix features and the optimal $\theta$:
\begin{align*}
\theta_{opt} = g(X_g) \ ,
\end{align*}
where $g$ is the mapping constructed through MLP. Let $\theta_{auto}$ denote the predicted value of $\theta$ by AutoAMG($\theta$), and $\theta_{opt}$ denote the optimal value of $\theta$. We use MSE (Mean Squared Error)~\cite{bickel2015mathematical} function as the loss function, then the loss is defined as 
\begin{equation}
\begin{split}
Loss &= MSE(\theta_{opt},\ \theta_{auto}) \\
     &= \frac{1}{M} \sum_{i=1}^{M} (\theta_{opt,i} - \theta_{auto,i})^2	
\end{split} \label{eq-mse}
\end{equation}
where $M$ is the batch size, $\theta_{opt,i}$ is the optimal $\theta$ value of the $i$-th matrix in the batch and $\theta_{auto,i}$ is the predicted $\theta$ value of the $i$-th matrix in the batch. The program of GCIN and optimization are implemented by PyTorch Geometric~\cite{2019pyg}.

\section{Numerical experiments}
\label{sec4}

\begin{figure}[h] 
	\centering
	\includegraphics[scale=1.8]{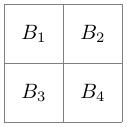}
	\caption{An example of $bx = by = 2$ blocks ($B_i, i = {1,2,3,4}$) with equal size. Diffusion coefficent $\kappa$ is the same in each block, when $B_i\neq B_j$, $\kappa_i \neq \kappa_j$.}
	\label{fig-block}
\end{figure}

We validated the effectiveness of AutoAMG($\theta$) based on two typical types of problems: the diffusion equation (\ref{eq-1}) with random coefficients, and the three-dimensional radiation diffusion equations~\cite{xiaowen2009algebraic} from inertial confinement fusion (as described in Section \ref{sec4-rde}). For Eq. (\ref{eq-1}), considering the 2D and 3D cases, the domain is $[0,1]^d (d=2,3)$, and the diffusion coefficients are
\begin{align}
\kappa = \begin{bmatrix}
10^{M r_0} & 0 \\
0 & 10^{M r_1}
\end{bmatrix} , \quad 
\kappa = \begin{bmatrix}
10^{M r_0} & 0 & 0\\
0 & 10^{M r_1} & 0 \\
0 & 0 & 10^{M r_2}
\end{bmatrix},  \label{eq-rand}
\end{align} 
where $r_0$, $r_1$, and $r_2$ are random numbers in the interval $(0,1)$, and $M \in \mathbb{N}_+$ is the parameter that influences the multiscale property of the matrix. A larger value of $M$ generally leads to a more pronounced multiscale property within the generated matrix. The computational domain is uniformly divided into blocks or subdomains with equal size, as shown in Figure \ref{fig-block}. While the diffusion coefficient $\kappa$ remains consistent within each block, it differs between different blocks. Therefore, even with identical mesh size and block count, different random seeds can generate distinct matrices.

When discretizing Eq. (\ref{eq-1}), matrices with varied properties and sizes can be generated by selecting different random number seeds $Seed$, mesh sizes $nx,ny,nz$ in each axis direction; block counts $bx,by,bz$ in each axis direction, and the parameter $M$. The matrix data is obtained from the following three equations: 
\begin{itemize}
	\item 2D diffusion equations: $nx=ny \in (50,100)$, $bx=by \in (10,20)$, $M=5$, and random seed $Seed$ is equal to the index of the matrix.
	\item 3D diffusion equations: $nx=ny=nz \in (30,40)$, $bx=by=bz \in (10,20)$, $M=5$, and random seed $Seed$ is equal to the index of the matrix.
	\item 3D radiation diffusion equations: 10 matrices from 3D radiation diffusion equations, with each matrix having approximately $6.29*10^6$ rows.
\end{itemize}

The optimal $\theta$ for each matrix is determined through grid search. We calculate the number of iterations by considering values of $\theta$ in increments of 0.01 within the range of $[0.01,0.99]$. The optimal $\theta$ is chosen as the one that results in the minimum number of iterations. The linear equations are solved using the JXPAMG software~\cite{xu2022jxpamg}, utilizing the GMRES algorithm with the AMG preconditioner. The coarsening algorithm in AMG is PMIS. We set an upper limit of 500 iterations, and the stopping criterion is that the relative residual is less than $10^{-8}$.

\begin{figure} 
	\centering 
	\includegraphics[scale=0.6]{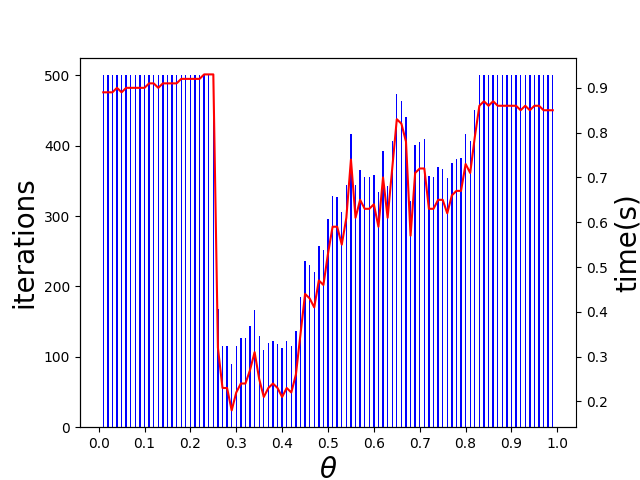}
	\caption{The matrix is from 2D diffusion equation with 9409 rows. Histogram is the number of iterations with left $y$-axis, and red line is the elapsed time with right $y$-axis}
	\label{fig-itertime}
\end{figure}

\begin{remark}
	The number of iterations is selected as the performance metric. While considering the operator complexity of AMG is closely related to the value of $\theta$, the elapsed time may seem like a preferable alternative. However, after plotting the number of iterations and time in the same picture (Figure \ref{fig-itertime}), it becomes evident that their trends are quite similar. Furthermore, given the matrix sizes in our experiments, some elapsed times are too short for precise measurement and are susceptible to the runtime environment. In contrast, the number of iterations remains unaffected by the environment. Hence, we have decided to utilize the number of iterations as our primary metric.
\end{remark}

\subsection{2D diffusion equations}
\label{sec4-2d}

The training and test sets consist of 80 and 20 matrices respectively. The mesh size $nx=ny \in (50,100)$ and the number of blocks $bx=by \in (10,20)$ are both random values. The results of the test set are shown in Table \ref{tab-exp2}. The first column "nrow" is the average number of rows of matrices in the test set; "iter" is the average number of iterations; "time" is the average time used to solve a linear equation in the test set. In the column of AutoAMG($\theta$), the "iter" and "time" correspond to the average number of iterations and computation time based on the $\theta$ predicted by AutoAMG($\theta$). The column "speedup" is the time of $\theta=0.25$ divided by the time of AutoAMG($\theta$).

\begin{table}[h] 
	\centering
	\caption{Test results of 2D diffusion equations}
	\label{tab-exp2}
	\begin{tabular}{cccccccc}
		\toprule
		\multirow{2}{*}{nrow} & \multicolumn{2}{c}{optimal $\theta$} & \multicolumn{2}{c}{$\theta = 0.25$} & \multicolumn{2}{c}{AutoAMG($\theta$)} & \multirow{2}{*}{speedup} \\
		\cmidrule(lr){2-3} \cmidrule(lr){4-5} \cmidrule(lr){6-7} 
	    & iter & time(s) & iter & time(s) & iter & time(s) &  \\
		\midrule
		5658.50 & 185.25 & 0.15 & 496.20 & 0.38 & 257.30 &  0.21 & 1.81 \\  
		\bottomrule
	\end{tabular}
\end{table}

Our objective is to assess the solving efficiency of the $\theta$ predicted by AutoAMG($\theta$) in comparison to the default $\theta$. Given that the default value of $\theta$ for 2D problems is 0.25, Table \ref{tab-exp2} presents the number of iterations and time corresponding to $\theta=0.25$. Despite the improved solving efficiency achieved by AutoAMG($\theta$), a noticeable gap remains between the attained performance and the optimal one.

\subsection{3D diffusion equations}
\label{sec4-3d}

The training and test sets consist of 80 and 20 matrices respectively. The mesh size $nx=ny=nz \in (30,40)$ and number of blocks $bx=by=bz \in (10,20)$ are random values. The results of the test set are shown in Table \ref{tab-exp3}, and the notations used are similar to those in Table \ref{tab-exp2}. In 3D equations, The number of iterations and time tuned by AutoAMG($\theta$) are close to the optimal ones, which is a significant improvement over the default value $\theta=0.5$.

\begin{table}[h] 
	\centering
	\caption{Test results of 3D diffusion equations}
	\label{tab-exp3}
	\begin{tabular}{cccccccc}
	\toprule
	\multirow{2}{*}{nrow} & \multicolumn{2}{c}{optimal $\theta$} & \multicolumn{2}{c}{$\theta = 0.5$} & \multicolumn{2}{c}{AutoAMG($\theta$)} & \multirow{2}{*}{speedup} \\
	\cmidrule(lr){2-3} \cmidrule(lr){4-5} \cmidrule(lr){6-7} 
	& iter & time(s) & iter & time(s) & iter & time(s) &  \\
	\midrule
	40514.90 & 34.00 & 0.29 & 233.20 & 1.52 & 42.75 & 0.34 & 4.47 \\  
	\bottomrule
	\end{tabular}
\end{table}

\subsection{Mixed 2D and 3D diffusion equations}
\label{sec4-2d3d}

A more common scenario arises when the origin of a matrix is unknown, making it challenging to determine whether it was discretized from a 2D or 3D problem. In such cases, AutoAMG($\theta$) is required to process the input matrix without additional information. Matrices from 2D and 3D diffusion equations are combined to make up the training and test sets, comprising 160 and 40 matrices respectively. To ensure a balanced distribution of matrix data, half of the data originates from 2D problems and the remaining half from 3D problems, both in the training and test sets. Since the dimension is unknown, we calculate the average number of iterations and computation time for all matrices in the test set at $\theta=0.25$ and $\theta=0.5$, as displayed in Table \ref{tab-mixed}.

\begin{table}[h] 
	\centering
	\caption{Test results of the mixed problems}
	\label{tab-mixed}
	\begin{tabular}{cccccccccc}
		\toprule
		\multicolumn{2}{c}{optimal $\theta$} & \multicolumn{2}{c}{$\theta = 0.25$} & \multicolumn{2}{c}{$\theta = 0.5$} & \multicolumn{2}{c}{AutoAMG($\theta$)} & \multicolumn{2}{c}{speedup} \\
		\cmidrule(lr){1-2} \cmidrule(lr){3-4} \cmidrule(lr){5-6} \cmidrule(lr){7-8} \cmidrule(lr){9-10}
		iter & time(s) & iter & time(s) & iter & time(s) & iter & time(s) & 0.25 & 0.5 \\
		\midrule
		109.63 & 0.22 & 273.83 & 0.39 & 291.00& 0.90& 179.88 & 0.29 & 1.34 & 3.10 \\
		\bottomrule
	\end{tabular}
\end{table}

From Table \ref{tab-mixed}, it is evident that the predicted $\theta$ by AutoAMG yields higher solving efficiency compared to default values of $\theta=0.25$ and $\theta=0.5$. However, training with mixed matrices results in a less robust model. The speedup over $\theta = 0.25$ and $\theta = 0.5$ is 1.34 and 3.10, whereas the speedup in Table \ref{tab-exp2} and \ref{tab-exp3} are 1.81 and 4.47. Consequently, it is advisable to train the model using matrices from the same dimension.

\subsection{3D radiation diffusion equations}
\label{sec4-rde}

The matrices employed in the previous sections originate from diffusion equations, containing fewer than $5 \times 10^4$ rows. To ascertain the generalizability of AutoAMG($\theta$), we employ all matrices from Section \ref{sec4-3d} for training and 10 matrices discretized from 3D radiation diffusion equations (with approximately $6.29 \times 10^6$ rows) for testing. Experimental results confirm that AutoAMG($\theta$) can be trained on smaller matrices and subsequently applied to larger matrices.

The formulas of 3D radiation diffusion equations \cite{xu2017algebraic,huang2022alphasetup} are
\begin{equation}
\begin{split}
c_{vr} \frac{\partial T_r}{\partial t} -\frac{1}{\rho} \nabla \cdot (K_r \nabla T_r) &= \omega_{er} (T_e- T_r) \ ,\\
c_{ve} \frac{\partial T_e}{\partial t} -\frac{1}{\rho} \nabla \cdot (K_e \nabla T_e) &= \omega_{ei} (T_i- T_e) + \omega_{er} (T_r - T_e) \ , \\
c_{vi} \frac{\partial T_i}{\partial t} -\frac{1}{\rho} \nabla \cdot (K_i \nabla T_i) &= \omega_{ei} (T_e- T_i)  \ ,
\end{split}
\label{eq-rde}
\end{equation}
where $\rho$ is the density; $T_r, T_e, T_i$ are the temperatures of photons, electrons, and ions, respectively; $c_{vr}, c_{ve}, c_{vi}$ are the specific heat at constant volume of photons, electrons, and ions, respectively; $K_r = f_r(\rho,T_r)$, $K_e = f_e(\rho,T_e)$ and $K_i = f_i(\rho,T_i)$ ($f_r, f_e, f_i$ are functions) are diffusion coefficients; and $\omega_{ei}$ and $\omega_{er}$ are the respective energy exchange coefficients between electrons and ions, and electrons and photons. Eq. (\ref{eq-rde}) is a nonlinear partial differential equation. It is discretized in time by the backward Euler method, then the nonlinear problem is transformed into a linear problem by the coagulation coefficient method, and the linear problem is discretized by the finite volume method. The sparse pattern of the discretized matrix is
\begin{align}
A = \begin{bmatrix}
A_R & D_{RE} & 0 \\
D_{ER} & A_E & D_{EI} \\
0 & D_{IE} & A_I 
\end{bmatrix}\ .
\label{eq-mat}
\end{align}
The block matrices $A_R$, $A_E$, $A_I$ in Eq. (\ref{eq-mat}) have the same sparse pattern, and the block matrices $D_{RE}$, $D_{EI}$ are diagonal matrices.

The training set consists of 100 matrices from Section \ref{sec4-3d} (80 matrices from the training set and 20 matrices from the test set), while the test set includes 10 matrices from Eq. (\ref{eq-rde}). The results are shown in Table \ref{tab-rde}. The number of iterations and computation time based on the $\theta$ predicted by AutoAMG($\theta$) are close to optimal ones, which is a substantial improvement compared to the default $\theta=0.5$. Moreover, in contrast with the speedup shown in Table \ref{tab-exp3}, AutoAMG($\theta$) demonstrates the capability to achieve even greater speedup. Such results illustrate the benefit of tuning $\theta$ in practical problems.

\begin{table}[h] 
	\centering
	\caption{Test results of 3D radiation diffusion equations}
	\label{tab-rde}
	\begin{tabular}{cccccccc}
		\toprule
		\multirow{2}{*}{nrow} & \multicolumn{2}{c}{optimal $\theta$} & \multicolumn{2}{c}{$\theta = 0.5$} & \multicolumn{2}{c}{AutoAMG($\theta$)} & \multirow{2}{*}{speedup} \\
		\cmidrule(lr){2-3} \cmidrule(lr){4-5} \cmidrule(lr){6-7} 
		& iter & time(s) & iter & time(s) & iter & time(s) &  \\
		\midrule
		6291456.00 & 31.50 & 31.52 & 484.20 & 399.00 & 35.40 & 34.27 & 11.63 \\  
		\bottomrule
	\end{tabular}
\end{table} 

\begin{figure}
	\centering
	\includegraphics[scale=0.5]{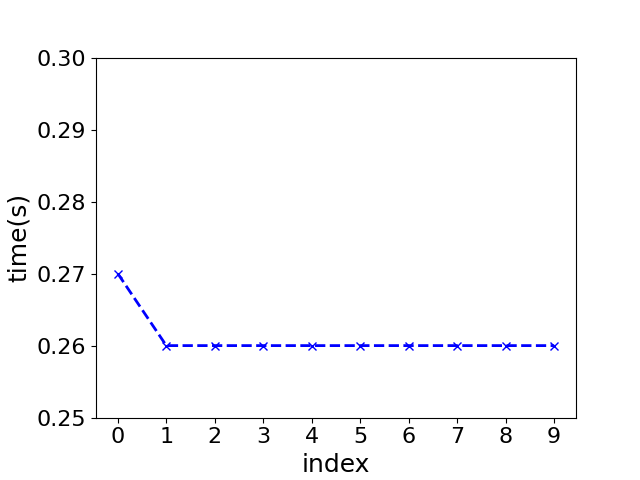}
	\caption{The inference time of each matrix in the test set}
	\label{fig-time}
\end{figure}

In terms of the overhead induced by AutoAMG($\theta$), we measure the inference time of each matrix in the test set, and the results are shown in Figure \ref{fig-time}. The $x$ axis in the figure is the index of the matrix, and the $y$ axis is the inference time. Note that the average inference time is 0.26 s, which is negligible compared to the average solving time of 34.27 s in Table \ref{tab-rde}. In fact, according to Eq. (\ref{eq-gcin1}), the message-passing process of GCIN can be effectively realized through the SpMV operation, hence it's conceivable that the overhead of GCIN would be inconsequential.

\section{Summary}
\label{sec5}

In this paper, we propose AutoAMG($\theta$), an auto-tuning method designed to adaptively adjust the strong threshold $\theta$ in the AMG algorithm for matrices from different problems. The effectiveness of this method is verified through a variety of numerical experiments.

An innovative contribution of this paper is the introduction of the GCIN algorithm for extracting matrix features. In diffusion problems, when compared to default $\theta$, the AutoAMG($\theta$) method based on GCIN demonstrates a speedup by a factor of 1.81 in 2D diffusion problems and 4.47 in 3D diffusion problems. Furthermore, AutoAMG($\theta$) displays versatility by effectively handling matrices from both 2D and 3D problems. Although it shows superior efficiency compared to default values, the speedup is only 1.34 in 2D problems and 3.10 in 3D problems.
 
Notably, in 3D radiation diffusion problems, AutoAMG($\theta$) effectively tunes the number of iterations and time that are close to the optimal results, achieving an impressive acceleration by a factor of 11.63 over the default $\theta=0.5$. The experiments reveal that AutoAMG($\theta$) generalizes well to new large matrices after training on small matrices.

Our future research will continue to focus on AMG algorithm optimization, using GNN to optimize the smoothing, interpolation, restriction, and other operators in AMG.

%%%% Acknowledgments %%%%%%%%
\section*{Acknowledgments} 
This work is financially supported by the National Natural Science Foundation of China (62032023) .

%Bibliography
\bibliographystyle{unsrt}  
\bibliography{autoamg}

\end{document}